\documentclass{article}

     \PassOptionsToPackage{numbers, compress}{natbib}

\usepackage{neurips_ml4co_2021}

\usepackage[utf8]{inputenc} 
\usepackage[T1]{fontenc}    
\usepackage{hyperref}       
\usepackage{url}            
\usepackage{booktabs}       
\usepackage{amsfonts}       
\usepackage{nicefrac}       
\usepackage{microtype}      
\usepackage{xcolor}         

\usepackage{amsmath,amssymb}
\usepackage{csquotes}
\usepackage{mathtools}
\usepackage{graphicx}
\usepackage{threeparttable}
\usepackage{multirow, multicol}
\usepackage{booktabs}

\graphicspath{{figures/}}

\newtheorem{thm}{Theorem}
\newtheorem{obs}[thm]{Observation}

\title{Efficient Primal Heuristics for Mixed-Integer Linear Programs}

\author{
	Akang Wang, Linxin Yang, Sha Lai, Xiaodong Luo \\
	Shenzhen Research Institute of Big Data\\
	Shenzhen, China \\
	\texttt{wangakang@sribd.cn, yanglinxin, laisha, xiaodongluo@cuhk.edu.cn} \\
	\AND
	Xiang Zhou, Haohan Huang, Shengcheng Shao, Yuanming Zhu, Dong Zhang, Tao Quan \\
	Huawei GTS, China \\
	\texttt{zhouxiang60, huanghaohan, shaoshengcheng, zhuyuanming5@huawei.com} \\
	\texttt{zhangdong48, quantao@huawei.com} 
}

\begin{document}

\maketitle

\begin{abstract}
This paper is a short report about our work for the primal task in the Machine Learning for Combinatorial Optimization NeurIPS 2021 Competition. 
For each dataset of our interest in the competition, we propose customized primal heuristic methods to efficiently identify high-quality feasible solutions.
The computational studies demonstrate the superiority of our proposed approaches over the competitors'.
\end{abstract}

\section{Introduction}
Many combinatorial optimization problems can be modeled and solved via \textit{Mixed-integer linear programming}~(MILP). 
As shown by~(\ref{eq:MILP}), the MILP problem seeks for an assignment to continuous and discrete decision variables ($x$ and $y$, respectively) that are subject to linear constraints, with the goal of minimizing a linear objective.
A wide range of important applications in transportation, manufacturing, and many other industrial domains can be formulated as MILP models.
However, solving those models to guaranteed optimality is generally $\mathcal{NP}$-hard.
\begin{equation}
	\begin{aligned}
		& \underset{x, y}{ \text{min} }  && c^\top x + d^\top y      \\
		& \text{ s.t. }                  && A x + By \leq h   \\
		&                                && x \in \mathbb{R}^n, y \in \mathbb{Z}^m  
	\end{aligned} \label{eq:MILP}
\end{equation}

Various general-purpose primal heuristics have been developed to produce feasible solutions to MILP problems.
Among them, the simplest one is \textit{rounding heuristic}, which rounds each fractional element in a solution to its nearest integer.
Although the rounding method is computationally cheap, it often fails to identify feasible points.
A more sophisticated idea is called \textit{feasibility pump}~\cite{fischetti2005feasibility, bertacco2007feasibility}.
Briefly speaking, feasibility pump iteratively calls a rounding step and a projection step to ensure that the returned solutions satisfy integrality constraints and linear constraints, respectively.
Usually, the aforementioned heuristics could not produce solutions of high quality and thus are often employed to construct the very first feasible solutions, while RINS~\cite{danna2005exploring}, RENS~\cite{berthold2014rens}, and other computationally expensive primal heuristics are used for improving the incumbent solutions.
These improvement heuristics search for new incumbents within a neighborhood of some given feasible solutions.

In this work, we tackle a few problem-specific MILP formulations of our interest by utilizing classic primal heuristics but in a more judicious manner, by exploiting the problem structures. 
We provide problem definitions in Section~\ref{sec:problem_definition} and our solution approaches in Section~\ref{sec:solution_approaches}.
Computational results are summarized in Section~\ref{sec:results}.
Finally, we conclude our work in Section~\ref{sec:conclusion}.

\section{Problem definitions}  \label{sec:problem_definition}
\begin{figure}[!h]
	\centering
	\includegraphics[scale=0.75]{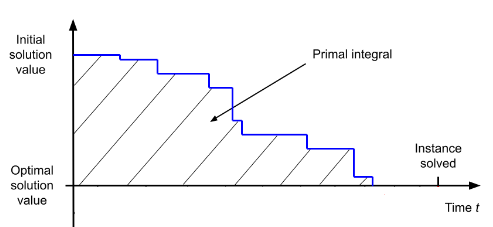}
	\caption{Primal bounds evolution versus solving time (adapted from \cite{ml4co})}
	\label{fig:primal_integral}
\end{figure}

The Machine Learning for Combinatorial Optimization (ML4CO) NeurIPS 2021 Competition~\cite{ml4co} focuses on the design of application-specific algorithms for solving MILP models.
Specifically, the primal task aims to identify high-quality incumbent solutions as fast as possible.
For an MILP instance of the form~(\ref{eq:MILP}), as a primal algorithm proceeds, new improved solutions will be identified and thus primal bounds are updated, as shown in Figure~\ref{fig:primal_integral}. 
The primal integral (the area of the shaded region) is then used as the metric for evaluating different primal algorithms. 
Clearly, an algorithm performing well will exhibit a small primal integral.

There are three problem benchmarks in the ML4CO competition: \textit{item placement}, \textit{load balancing}, and \textit{anonymous}. 
The first two problem benchmarks are inspired by real-life applications of large-scale systems at Google, while the third benchmark is an anonymous problem inspired by a real-world, large-scale industrial application. 
There are $10, 000$ instances for the first two benchmarks, respectively, and $118$ instances in the anonymous dataset, as training data. 
The competition encourages both machine learning (ML) based algorithms and non-ML algorithms.
\subsection{Item placement}
This item placement problem deals with spreading items across multiple containers utilizing them evenly. 
Let $I$ denote the set of items and $J$ denote the set of containers. 
Let a binary variable $x_{ij}$ be $1$ if item $i$ is placed in container $j$ and $0$ otherwise.
Each item will be placed in exactly a single container, as shown by constraints~(\ref{eq:ip:assignment}). 
Let $K$ represent the set of dimensions.
For dimension $k \in K$ of container $j \in J$, knapsack constraints~(\ref{eq:ip:knapsack}) represent some physical considerations while (\ref{eq:ip:define_y}) and~(\ref{eq:ip:define_z}) properly account for the placement unevenness, which is penalized in the objective~(\ref{eq:ip:obj}).  
The goal is to identify an assignment for $x, y, z$ such that constraints~(\ref{eq:ip:assignment})~--~(\ref{eq:ip:positive}) are satisfied and the objective~(\ref{eq:ip:obj}) is minimized.
	\begin{align}
		& \underset{x, y, z}{\text{min}} && \sum_{j \in J}^{} \sum_{k \in K}^{} \alpha_{k} y_{jk} + \sum_{k \in K}^{} \beta_k z_{k} \label{eq:ip:obj} \\
		& \; \text{s.t.}  && \sum_{j \in J}^{} x_{ij} = 1 &&& \forall i \in I \label{eq:ip:assignment} \\
		&              && \sum_{i \in I }^{} a_{ik} x_{ij} \leq b_{k}    &&& \forall j \in J, \forall k \in K  \label{eq:ip:knapsack}  \\
		&              && \sum_{i \in I}^{} d_{ik} x_{ij} + y_{jk}  \geq 1   &&& \forall j \in J, \forall k \in K \label{eq:ip:define_y} \\
		&              && y_{jk} \leq z_{k}    &&& \forall j \in J, \forall k \in K  \label{eq:ip:define_z} \\
		&              && x_{ij} \in \left\{0, 1\right\} &&& \forall i \in I, \forall j \in J \\
		&              && y_{jk} \geq 0         &&& \forall j \in J, \forall k \in K  \label{eq:ip:positive}
	\end{align}

\subsection{Load balancing}
This problem deals with apportioning workloads across as few machines as possible. 
Let $M$ and $N$ denote the set of tasks and the set of machines, respectively.
For task $i \in M$, only a subset of machines, denoted by $N^i \subseteq N$, are accessible. 
Let a binary variable $y_j$ be $1$ if machine $j$ is used and $0$ otherwise.
Let $x_{ij}$ denote the amount of workload from task $i$ to machine $j$, as defined in constraints~(\ref{eq:lb:define_x}).
Constraints~(\ref{eq:lb:capacity}) enforce the capacity requirement for each machine.
The apportionment is required to be robust to any one machine's failure, as indicated by constraints~(\ref{eq:lb:robust}). 
The goal is to identify an assignment for $x, y$ such that constraints~(\ref{eq:lb:define_x})~--~(\ref{eq:lb:positive}) are satisfied and the objective~(\ref{eq:lb:obj}) is minimized.
	\begin{align}
		& \underset{x, y}{\text{min}} && \sum_{j \in N}^{} y_j  \label{eq:lb:obj} \\
		& \; \text{s.t.}  &&  x_{ij} \leq a_i y_j &&& \forall i \in M, \forall j \in N^i  \label{eq:lb:define_x}\\
		&              && \sum_{i \in M: j \in N^i}^{} x_{ij} \leq b_j    &&& \forall j \in N   \label{eq:lb:capacity}  \\
		&              && \sum_{j \in N^i \setminus \left\{j^\prime\right\}}^{}  x_{ij} \geq a_i   &&& \forall i \in M, \forall j^\prime \in N^i  \label{eq:lb:robust}  \\
		&              && y_j \in \left\{0, 1\right\}   &&& \forall j \in N \\
		&             && 0 \leq x_{ij} \leq b_{j}  &&& \forall i \in M, \forall j \in N^i  \label{eq:lb:positive}
	\end{align} 

\subsection{Anonymous}
The concrete MILP formulation corresponding to this benchmark is not provided in the ML4CO competition.

\section{Solution approaches} \label{sec:solution_approaches}

\subsection{Item placement}
We analyze $10, 000$ item placement instances and find out that 
\begin{itemize}
	\item $|I| = 105, |J| = 10$;
	\item $a_{ik}, d_{ik}$ values of $5$ items are significantly bigger than those of others. 
\end{itemize}
Based on this empirical finding, we place the $5$ big items into the first five containers respectively before any primal heuristics are applied. 

\paragraph{Meta-heuristics}
The item placement problem is simply a multi-dimensional multi-knapsack instance.
The knapsack problem is one of the most well-studied combinatorial optimization problems, and tremendous efforts have been devoted to the development of efficient heuristic algorithms, from which we adopt some ideas for item placement. 
In particular, we first apply a greedy method in which items are first sorted based on their sizes and then assigned to containers.
This will produce the very first feasible solution.
After that, a \textit{large neighborhood search} idea is considered as an improvement heuristic approach. 
We select one or two items respectively from two containers and swap them if this leads to a better incumbent.  

\paragraph{Math-heuristics} Math-heuristic methods are based on solving mathematical models. 
We consider a construction method and an improvement method.
The construction method consists of two steps:
(i)~first assign items to the first five containers by solving an assignment model; 	(ii)~then assign the remaining items to the last five containers by solving a sub-MIP.
The assignment model is defined as follows:
\begin{itemize}
	\item[(a)] combine small items to generate candidates that might be placed in the first five containers;
	\item[(b)]  assume all items that are placed at the remaining containers have \enquote{continuous} sizes (rather than \enquote{discrete});
	\item[(c)] add an SOS1 constraint for each small item (to ensure that each item will be placed in at most one container).
\end{itemize}

Again, the large neighborhood search method is utilized for solution improvement. 
Specifically, we properly choose two out of the last five containers, and then solve a sub-MIP to reassign items within those two containers optimally.

\begin{figure}[!h]
	\centering
	\includegraphics[scale=0.42]{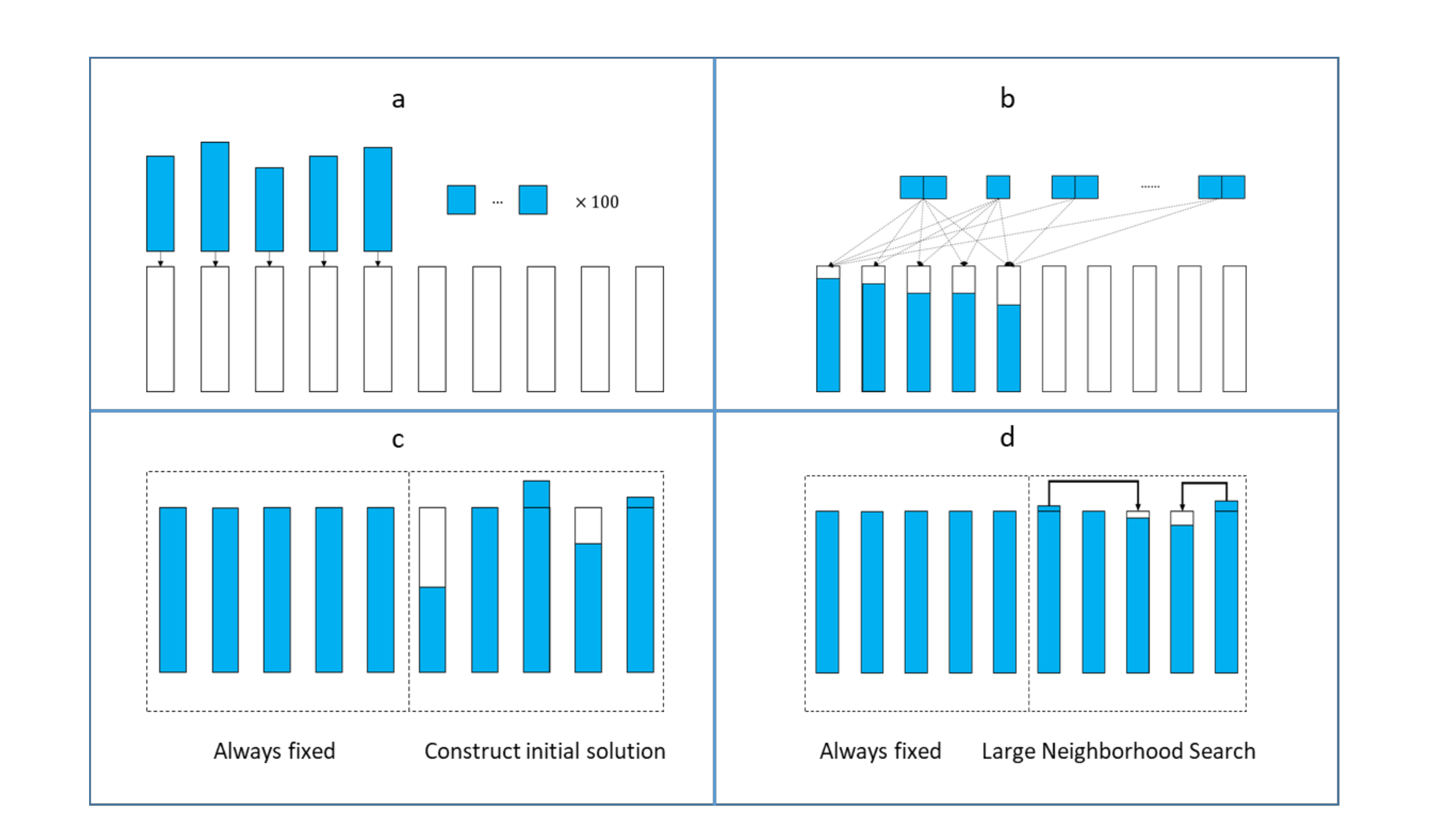}
	\caption{Math-heuristics for item placement}
	\label{fig:item_placement}
\end{figure}

A pictorial representation of the math-heuristic method is given in Figure~\ref{fig:item_placement}.

\subsection{Load balancing}
\begin{obs}  \label{obs:rounding_up}
	Rounding up a solution to the linear programming relaxation of model (\ref{eq:lb:obj})~--~(\ref{eq:lb:positive}) would produce a feasible solution to (\ref{eq:lb:obj})~--~(\ref{eq:lb:positive}).
\end{obs}
Observation~\ref{obs:rounding_up} is clearly valid since rounding up the fractional element $y_j$ with a non-zero value will not affect the satisfaction of linear constraints in the model but produce an integral solution.
To further exploit the possibility of rounding a fractional solution towards a new incumbent, in our implementation, we adaptively select the rounding threshold $\eta$ for $y_j$.
In another word, we choose $\eta$ in an adaptive manner and round up $y_j$ only if it exceeds $\eta$. 
The large $\eta$ values would lead to solutions with smaller objective values but might cause returned solutions to be infeasible, hence some care has been taken when choosing $\eta$.
Specifically, we first select a target objective value based on the current primal and dual bounds and then determine $\eta$ via quantile selection such that after rounding with $\eta$ the objective matches the predetermined value. 
The initial primal bound corresponds to the feasible solution returned by rounding up every fractional element $y_j$ with a non-zero value, while the initial dual bound is approximated by amplifying the root node lower bound with a certain factor. 
If the rounding step produces a new incumbent solution, we then update the primal bound; otherwise, we set the dual bound to the corresponding objective value. 
We iterate this process until the primal-dual gap falls below a certain threshold.

\begin{obs}  \label{obs:tightening}
	In model~(\ref{eq:lb:obj})~--~(\ref{eq:lb:positive}), constraints~(\ref{eq:lb:capacity}) can be tightened as follows:
	\begin{equation}
		\begin{aligned}
			\sum_{i \in M: j \in N^i}^{} x_{ij} \leq b_j y_j    &&& \forall j \in N. 
		\end{aligned}  \label{eq:lb:tightened_capacity}
	\end{equation}
\end{obs}
Due to constraints~(\ref{eq:lb:define_x}), $y_j = 0$ implies $x_{ij} = 0$, hence Observation~\ref{obs:tightening} is valid. 
Furthermore, across $10, 000$ instances, we also observe that $a_i < b_j$ for $i \in M, j \in N^i$. 
As a result, constraints~(\ref{eq:lb:define_x}) are dominated by~(\ref{eq:lb:tightened_capacity}) and thus eliminated from the model.
Now we call~(\ref{eq:lb:obj}), (\ref{eq:lb:robust})~--~(\ref{eq:lb:tightened_capacity}) as the \enquote{tightened model}. 

In our experience, the tightened model has a much better root node bound. 
Thus, the root node linear programming~(LP) solution is likely to be a good candidate for the rounding heuristics we discussed previously. 
In our implementation, we first apply the rounding heuristic method to the root LP solution to the original model and then to the optimal LP solution to tightened model.
To further improve the incumbent, we then adopt RINS. 
In particular, we define and solve a sub-MIP based on the tightened model, using its LP solution and the current incumbent as a guide to fix part of the binary variables. 

\subsection{Anonymous}
\begin{figure}[!h]
	\centering
	\includegraphics[scale=0.6]{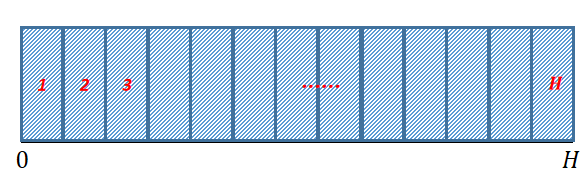}
	\caption{Planning horizon}
	\label{fig:planning_horizon}
\end{figure}

Though the concrete MILP formulation is not available, we can discover some pattern from MILP instances in the LP file format.
\begin{obs} \label{obs:anonymous}
	One can define a planning horizon and associate with each discrete variable a time period.
\end{obs}  
The details of Observation~\ref{obs:anonymous} can be deduced from the constraint hypergraph~\cite{rossi2006handbook} in which every node represents a discrete variable and every edge joins a pair of variables if they occur together in a constraint.
Let $H$ denote the planning horizon and $h \in H$ denote a time period, as shown in Figure~\ref{fig:planning_horizon}.

\begin{figure}[!h]
	\centering
	\includegraphics[scale=0.6]{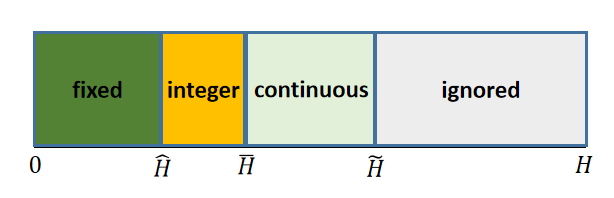}
	\caption{The rolling-horizon method}
	\label{fig:rolling_horizon}
\end{figure}

Based on Observation~\ref{obs:anonymous}, we propose our first primal heuristic called a \enquote{rolling-horizon} method (Figure~\ref{fig:rolling_horizon}), which consists of the following steps:
\begin{itemize}
	\item[(i)] ignore constraints involving discrete variables with their $h$ values greater than $\widetilde{H}$; 
	\item[(ii)] relax the integrality constraints on variables with their $h$ values greater than $\overline{H}$;
	\item[(iii)] fix those discrete variables with their $h$ values less than $\widehat{H}$ at the optimal solution from a previous run;
	\item[(iv)] solve the sub-MIP;
	\item[(v)] increase $\widehat{H}, \overline{H}$ and $\widetilde{H}$ and then iterate steps (i) - (iv) until $\widehat{H} = H$. 
\end{itemize}
Parameters $\widehat{H}, \overline{H}$ and $\widetilde{H}$ are selected adaptively to balance feasibility, optimality, and the computational cost of solving sub-MIPs. 
We remark that the rolling-horizon method might be generally applicable to optimization models with discrete variables spanning consecutive time periods, for example, inventory management models. 

The rolling-horizon method is often quite time-consuming, hence we first utilize some computationally favorable methods to generate the very first solutions.
In particular, we implement a variant of feasibility pump that combines a rounding step and a projection step, using weak and strong perturbation to escape cycle or stalling.
After obtaining the first incumbent, we then call RENS.
Specifically, the RENS model is a sub-MIP defined by fixing those discrete variables with their $h$ values less than $0.9H$ at the incumbent.

\section{Computational Results}  \label{sec:results}
We implement our code in Python except that the meta-heuristic algorithm is implemented in C++.
The evaluation is conducted via Ecole~0.7.3~\cite{prouvost2020ecole}, using the PySCIPOpt~3.3.0~\cite{maher2016pyscipopt} interface to call SCIP~7.0.3~\cite{gamrath2020scip} as the LP and MILP solver. 
The test set includes $100$ instances for both item placement and load balancing, $20$ anonymous instances (run each instance with $5$ different seeds).    
The time limit is $300$ seconds for each instance and the average primal integral is reported as the performance metric. 

We adopt the final results from~\cite{ml4co} and provide the synopsized version in Table~\ref{table:results}.
The columns \enquote{Item placement}, \enquote{Load balancing}, and \enquote{Anonymous} report the average primal integral values across $100$ runs, respectively.
Compared with the second places ($2^{nd}$ in Table~\ref{table:results}) in item placement, load balancing and anonymous benchmarks, our proposed approaches produce $67.5\%$, $51.8\%$, and $7.7\%$ smaller primal integral values, respectively.
Clearly, our customized approaches significantly outperform those from other teams in the competition.
\begin{table}[htbp] 
	\caption{Computational results in the ML4CO competition}
	\label{table:results}
	\centering
	\renewcommand{\arraystretch}{1}
	\setlength\tabcolsep{3mm}
	\begin{tabular}{cccc}
		\toprule
		\multirow{2}{*}{Team}  &	\multicolumn{3}{c}{Primal integral}\\
		\cmidrule(l){2-4} 
		& Item placement & Load balancing & Anonymous\\
		\midrule
		Our work      & 358.43  & 1145.31 & 35220755.00 \\
		$2^{nd}$      & 1102.19 & 2374.55 & 38163724.79 \\
		\bottomrule
	\end{tabular}
\end{table}

\section{Conclusions}  \label{sec:conclusion}
In this work, we propose specialized primal algorithms for tackling item placement, load balancing and anonymous benchmarks in the ML4CO competition.
Our approaches are mainly based on classic primal heuristics in the MILP literature. 
However, we exploit the problem structures and employ those heuristics in a much more sophisticated manner.
As a result, our proposed algorithms perform significantly better than the other approaches in the ML4CO competition.

\bibliographystyle{plain}
\bibliography{main}

\begin{thebibliography}{1}

\bibitem{ml4co}
\url{https://www.ecole.ai/2021/ml4co-competition/}, accessed on January 22,
  2022.

\bibitem{bertacco2007feasibility}
Livio Bertacco, Matteo Fischetti, and Andrea Lodi.
\newblock A feasibility pump heuristic for general mixed-integer problems.
\newblock {\em Discrete Optimization}, 4(1):63--76, 2007.

\bibitem{berthold2014rens}
Timo Berthold.
\newblock Rens.
\newblock {\em Mathematical Programming Computation}, 6(1):33--54, 2014.

\bibitem{danna2005exploring}
Emilie Danna, Edward Rothberg, and Claude Le~Pape.
\newblock Exploring relaxation induced neighborhoods to improve mip solutions.
\newblock {\em Mathematical Programming}, 102(1):71--90, 2005.

\bibitem{fischetti2005feasibility}
Matteo Fischetti, Fred Glover, and Andrea Lodi.
\newblock The feasibility pump.
\newblock {\em Mathematical Programming}, 104(1):91--104, 2005.

\bibitem{gamrath2020scip}
Gerald Gamrath, Daniel Anderson, Ksenia Bestuzheva, Wei-Kun Chen, Leon Eifler,
  Maxime Gasse, Patrick Gemander, Ambros Gleixner, Leona Gottwald, Katrin
  Halbig, et~al.
\newblock The scip optimization suite 7.0.
\newblock 2020.

\bibitem{maher2016pyscipopt}
Stephen Maher, Matthias Miltenberger, Joao~Pedro Pedroso, Daniel Rehfeldt,
  Robert Schwarz, and Felipe Serrano.
\newblock Pyscipopt: Mathematical programming in python with the scip
  optimization suite.
\newblock In {\em International Congress on Mathematical Software}, pages
  301--307. Springer, 2016.

\bibitem{prouvost2020ecole}
Antoine Prouvost, Justin Dumouchelle, Lara Scavuzzo, Maxime Gasse, Didier
  Ch{\'e}telat, and Andrea Lodi.
\newblock Ecole: A gym-like library for machine learning in combinatorial
  optimization solvers.
\newblock {\em arXiv preprint arXiv:2011.06069}, 2020.

\bibitem{rossi2006handbook}
Francesca Rossi, Peter Van~Beek, and Toby Walsh.
\newblock {\em Handbook of constraint programming}.
\newblock Elsevier, 2006.

\end{thebibliography}

\end{document}